\documentclass{amsart}
\usepackage{amssymb}
\usepackage{graphics}
\usepackage{amscd}
\newtheorem{thm}{Theorem}[section]
\newtheorem{cor}[thm]{Corollary}

\newtheorem{prop}[thm]{Proposition}

\theoremstyle{definition}
\newtheorem{defn}[thm]{Definition}

\theoremstyle{remark}
\newtheorem{rem}[thm]{Remark}
\newtheorem{ex}[thm]{Example}
\numberwithin{equation}{section}

\begin{document}

\title{Knot Floer homology and Seifert surfaces}%
\author{Andr\'as Juh\'asz}%
\address{Department of Mathematics, Princeton University, Princeton, NJ 08544, USA}%
\email{ajuhasz@math.princeton.edu}%

\thanks{Research partially supported by OTKA grant no. T49449}

\date{\today}%
\begin{abstract}
Let $K$ be a knot in $S^3$ of genus $g$ and let $n>0.$ We show that
if $\text{rk}\,\widehat{HFK}(K,g) < 2^{n+1}$ (where $\widehat{HFK}$
denotes knot Floer homology), in particular if $K$ is an alternating
knot such that the leading coefficient $a_g$ of its Alexander
polynomial satisfies $|a_g| <2^{n+1},$ then $K$ has at most $n$
pairwise disjoint non-isotopic genus $g$ Seifert surfaces. For $n=1$
this implies that $K$ has a unique minimal genus Seifert surface up
to isotopy.
\end{abstract}

\maketitle
\section{Introduction and preliminaries}

If $S_1$ and $S_2$ are Seifert surfaces of a knot $K \subset S^3$
then $S_1$ and $S_2$ are said to be equivalent if $S_1 \cap X(K)$
and $S_2 \cap X(K)$ are ambient isotopic in the knot exterior $X(K)
= S^3 \setminus N(K),$ where $N(K)$ is a regular neighborhood of
$K.$ In \cite{Kakimizu2} Kakimizu assigned a simplicial complex
$MS(K)$ to every knot $K$ in $S^3$ as follows.

\begin{defn} \label{defn:1}
$MS(K)$ is a simplicial complex whose vertices are the equivalence
classes of the minimal genus Seifert surfaces of $K.$ The
equivalence classes $\sigma_0, \dots, \sigma_n$ span an $n$-simplex
if and only if for each $0 \le i \le n$ there is a representative
$S_i$ of $\sigma_i$ such that the surfaces $S_0, \dots, S_n$ are
pairwise disjoint.
\end{defn}

In \cite{Thompson} it is shown that the complex $MS(K)$ is always
connected. I.e., if $S$ and $T$ are minimal genus Seifert surfaces
for a knot $K$ then there is a sequence $S=S_1, S_2, \dots, S_k = T$
of minimal genus Seifert surfaces such that $S_i \cap S_{i+1} =
\emptyset$ for $0 \le i \le k-1.$

The main goal of this short note is to show that for a genus $g$
knot $K$ and for $n>0$ the condition
$\text{rk}\,\widehat{HFK}(K,g)<2^{n+1}$ implies $\dim MS(K) < n,$
consequently for $n=1$ the knot $K$ has a unique Seifert surface up
to equivalence. This condition involves the use of knot Floer
homology introduced by Ozsv\'ath and Szab\'o in \cite{OSz3} and
independently by Rasmussen in \cite{Ras}. However, when $K$ is
alternating then this condition is equivalent to $|a_g|<2^{n+1},$
where $a_g$ is the leading coefficient of the Alexander polynomial
of $K.$ The alternating case is already a new result whose statement
doesn't involve knot Floer homology. On the other hand, the proof of
this particular case seems to need sutured Floer homology
techniques, which is a generalization of knot Floer homology that
was introduced by the author in \cite{sutured}. At the time of
writing this paper there are very few results in knot theory which
can only be proved using Floer homology methods.

The above statement does not hold for $n=0$ since every knot has at
least one minimal genus Seifert surface. However, it was shown in
\cite{fibred} and \cite{decomposition} that
$\text{rk}\,\widehat{HFK}(K,g)<2$ implies that the knot $K$ is
fibred.

To a knot $K$ in $S^3$ and every $j \in \mathbb{Z}$ knot Floer
homology assigns a graded abelian group $\widehat{HFK}(K,j)$ whose
Euler characteristic is the coefficient $a_j$ of the Alexander
polynomial $\Delta_K(t).$ In \cite{OSz9} it is shown that if $K$ is
alternating then $\widehat{HFK}(K,j)$ is non-zero in at most one
grading, thus $\text{rk}\, \widehat{HFK}(K,j) = |a_j|.$

Next we are going to review some necessary definitions and results
from the theory of sutured manifolds and sutured Floer homology.
Sutured manifolds were introduced by Gabai in \cite{Gabai}.

\begin{defn} \label{defn:2}
A \emph{sutured manifold} $(M,\gamma)$ is a compact oriented
3-manifold $M$ with boundary together with a set $\gamma \subset
\partial M$ of pairwise disjoint annuli $A(\gamma)$ and tori
$T(\gamma).$ Furthermore, the interior of each component of
$A(\gamma)$ contains a \emph{suture}, i.e., a homologically
nontrivial oriented simple closed curve. We denote the union of the
sutures by $s(\gamma).$

Finally every component of $R(\gamma)=\partial M \setminus
\text{Int}(\gamma)$ is oriented. Define $R_+(\gamma)$ (or
$R_-(\gamma)$) to be those components of $\partial M \setminus
\text{Int}(\gamma)$ whose normal vectors point out of (into) $M$.
The orientation on $R(\gamma)$ must be coherent with respect to
$s(\gamma),$ i.e., if $\delta$ is a component of $\partial
R(\gamma)$ and is given the boundary orientation, then $\delta$ must
represent the same homology class in $H_1(\gamma)$ as some suture.

A sutured manifold is called \emph{taut} if $R(\gamma)$ is
incompressible and Thurston norm minimizing in $H_2(M,\gamma).$
\end{defn}

The following definition was introduced in \cite{sutured}.

\begin{defn} \label{defn:3}
A sutured manifold $(M,\gamma)$ is called \emph{balanced} if M has
no closed components, $\chi(R_+(\gamma))=\chi(R_-(\gamma)),$ and the
map $\pi_0(A(\gamma)) \to \pi_0(\partial M)$ is surjective.
\end{defn}

\begin{ex} \label{ex:1} If $R$ is a Seifert surface of a knot $K$ in
$S^3$ then we can associate to it a balanced sutured manifold
$S^3(R) = (M,\gamma)$ such that $M = S^3 \setminus (R \times I)$ and
$\gamma = K \times I.$ Observe that $R_-(\gamma) = R \times \{0\}$
and $R_+(\gamma) = R \times \{1\}.$ Furthermore, $S^3(R)$ is taut if
and only if $R$ is of minimal genus.
\end{ex}

Sutured Floer homology is an invariant of balanced sutured manifolds
defined by the author in \cite{sutured}, and is a common
generalization of the invariants $\widehat{HF}$ and $\widehat{HFK}.$
It assigns an abelian group $SFH(M,\gamma)$ to each balanced sutured
manifold $(M,\gamma).$ The following theorem is a special case of
\cite[Theorem 1.5]{decomposition}.

\begin{thm} \label{thm:1}
Let $K$ be a genus $g$ knot in $S^3$ and suppose that $R$ is a
minimal genus Seifert surface for $K.$ Then
$$SFH(S^3(R)) \approx \widehat{HFK}(K,g).$$
\end{thm}

A sutured manifold $(M,\gamma)$ is called a product if it is
homeomorphic to $(\Sigma \times I, \partial \Sigma \times I),$ where
$\Sigma$ is an oriented surface with boundary. If $(M,\gamma)$  is a
product then $SFH(M,\gamma) \approx \mathbb{Z}.$ Let us recall
\cite[Theorem 1.4]{decomposition} and \cite[Theorem
9.3]{decomposition}.

\begin{thm} \label{thm:2}
If $(M,\gamma)$ is a taut balanced sutured manifold then
$SFH(M,\gamma) \ge \mathbb{Z}.$ Furthermore, if $(M,\gamma)$ is not
a product then $SFH(M,\gamma) \ge \mathbb{Z}^2.$
\end{thm}

\begin{defn} \label{defn:4}
Let $(M,\gamma)$ be a balanced sutured manifold. An oriented surface
$S \subset M$ is called a \emph{horizontal surface} if $S$ is open,
$\partial S=s(\gamma)$ in an oriented sense; moreover, $[S] =
[R_+(\gamma)]$ in $H_2(M,\gamma),$ and $\chi(S) =
\chi(R_+(\gamma)).$

A horizontal surface $S$ defines a horizontal decomposition
$$(M,\gamma) \rightsquigarrow^S (M_-,\gamma_-) \coprod
(M_+,\gamma_+)$$ as follows. Let $M_{\pm}$ be the union of the
components of $M \setminus \text{Int}(N(S))$ that intersect
$R_{\pm}(\gamma).$ Similarly, let $\gamma_{\pm}$ be the union of the
components of $\gamma \setminus \text{Int}(N(S))$ that intersect
$R_{\pm}(\gamma).$
\end{defn}

The following proposition is a special case of \cite[Proposition
8.6]{decomposition}.

\begin{prop} \label{prop:1}
Suppose that $(M,\gamma)$ is a taut balanced sutured manifold and
let $S$ be a horizontal surface in it. Then $$\text{rk}\,
SFH(M,\gamma) = \text{rk}\,SFH(M_-,\gamma_-) \cdot \text{rk} \,
SFH(M_+,\gamma_+).$$
\end{prop}

The following definition can be found for example in \cite{fibred}.

\begin{defn} \label{defn:5}
A balanced sutured manifold $(M,\gamma)$ is called
\emph{horizontally prime} if every horizontal surface $S$ in
$(M,\gamma)$ is isotopic to either $R_+(\gamma)$ or $R_-(\gamma)$
rel $\gamma.$
\end{defn}

\section{The results}

\begin{thm} \label{thm:3}
Let $(M,\gamma)$ be a taut balanced sutured manifold such that both
$R_+(\gamma)$ and $R_-(\gamma)$ are connected. Suppose that there is
a sequence of pairwise disjoint non-isotopic connected horizontal
surfaces $R_-(\gamma) = S_0, S_1, \dots, S_n = R_+(\gamma).$ Then
$$\text{rk}\,SFH(M,\gamma) \ge 2^n.$$
\end{thm}

\begin{proof}
We prove the theorem using induction on $n.$ If $n=1$ then
$(M,\gamma)$ is not a product since $R_-(\gamma)$ and $R_+(\gamma)$
are non-isotopic. Thus Theorem \ref{thm:2} implies that
$\text{rk}\,SFH(M,\gamma) \ge 2.$

Now suppose that the theorem is true for $n-1.$ Since each $S_k$ is
connected we can suppose without loss of generality that $S_1$
separates $S_i$ and $S_0$ for every $i \ge 2.$ Let $(M_-,\gamma_-)$
and $(M_+,\gamma_+)$ be the sutured manifolds obtained after
horizontally decomposing $(M,\gamma)$ along $S_1.$ Note that both
$(M_-,\gamma_-)$ and $(M_+,\gamma_+)$ are taut. Since $S_0$ and
$S_1$ are non-isotopic $(M_-,\gamma_-)$ is not a product so as
before $\text{rk}\, SFH(M_-,\gamma_-) \ge 2.$ Applying the induction
hypothesis to $(M_+,\gamma_+)$ and to the surfaces $R_-(\gamma_+),
S_2, \dots, S_n = R_+(\gamma_+)$ we get that
$\text{rk}\,SFH(M_+,\gamma_+) \ge 2^{n-1}.$ So using Proposition
\ref{prop:1} we see that $\text{rk}\,SFH(M,\gamma) \ge 2^n.$
\end{proof}

\begin{cor} \label{cor:1}
If $(M,\gamma)$ is a taut balanced sutured manifold and $\text{rk}\,
SFH(M,\gamma) < 4$ then $(M,\gamma)$ is horizontally prime. More
generally, if $n>0$ and $\text{rk}\,SFH(M,\gamma) < 2^{n+1}$ then
$(M,\gamma)$ can be cut into horizontally prime pieces by less than
$n$ horizontal decompositions.
\end{cor}

\begin{proof}
Suppose that $\text{rk}\,SFH(M,\gamma) < 2^{n+1}.$ If $(M,\gamma)$
is not horizontally prime then there is a surface $S_1$ in
$(M,\gamma)$ which is not isotopic to $R_{\pm}(\gamma).$ Decomposing
$(M,\gamma)$ along $S_1$ we get two sutured manifolds
$(M_-,\gamma_-)$ and $(M_+,\gamma_+).$ If they are not both
horizontally prime then repeat the the above process with a
non-prime piece and obtain a horizontal surface $S_2,$ etc. This
process has to end in less than $n$ steps according to Theorem
\ref{thm:3}.
\end{proof}

\begin{thm} \label{thm:4}
Let $K$ be a knot in $S^3$ of genus $g$ and let $n>0.$ If
$\text{rk}\,\widehat{HFK}(K,g) < 2^{n+1}$ then $K$ has at most $n$
pairwise disjoint non-isotopic genus $g$ Seifert surfaces, in other
words, $\dim MS(K) < n.$ If $n=1$ then $K$ has a unique Seifert
surface up to equivalence.
\end{thm}

\begin{proof}
Suppose that $R, S_1, \dots, S_n$ are pairwise disjoint non-isotopic
Seifert surfaces for $K.$ According to Theorem \ref{thm:1} we have
$\widehat{HFK}(K,g) \approx SFH(S^3(R)).$ Let $S^3(R) = (M,
\gamma).$ If $R_+(\gamma)$ and $R_-(\gamma)$ were isotopic then
$(M,\gamma)$ would be a product and $S_1$ and $R$ would be
equivalent. So the surfaces $R_-(\gamma)=S_0, S_1, \dots, S_n,
S_{n+1}=R_+(\gamma)$ satisfy the conditions of Theorem \ref{thm:3},
thus $\text{rk}\,SFH(S^3(R)) \ge 2^{n+1},$ a contradiction.

In particular, if $n=1$ then $\dim MS(K) =0.$ But according to
\cite{Thompson} the complex $MS(K)$ is connected, so it consists of
a single point.
\end{proof}

\begin{cor} \label{cor:2}
Suppose that $K$ is an alternating knot in $S^3$ of genus $g$ and
let $n>0.$ If the leading coefficient $a_g$ of its Alexander
polynomial satisfies $|a_g| < 2^{n+1}$ then $\dim MS(K) < n.$ If
$|a_g| < 4$ then $K$ has a unique Seifert surface up to equivalence.
\end{cor}

\begin{proof}
This follows from Theorem \ref{thm:4} and the fact that for
alternating knots $\text{rk}\,\widehat{HFK}(K,g) = |a_g|.$
\end{proof}

\begin{rem} \label{rem:1}
In \cite{Kakimizu} Kakimizu classified the minimal genus Seifert
surfaces of all the prime knots with at most $10$ crossings. The
$n=1$ case of Corollary \ref{cor:2} is sharp since the knot $7_4$ is
alternating, the leading coefficient of its Alexander polynomial is
$4,$ and has $2$ inequivalent minimal genus Seifert surfaces. On the
other hand, the Alexander polynomial of the alternating knot $9_2$
is also $4,$ but has a unique minimal genus Seifert surface up to
equivalence.

Also note that \cite[Thorem 1.7]{decomposition} implies that if the
leading coefficient $a_g$ of the Alexander polynomial of an
alternating knot $K$ satisfies $|a_g| < 4$ then the knot exterior
$X(K)$ admits a depth $\le 1$ taut foliation transversal to
$\partial X(K).$ Indeed, for alternating knots $g = g(K)$ and $|a_g|
= \text{rk}\, \widehat{HFK}(K,g) \neq 0,$ so the conditions of
\cite[Thorem 1.7]{decomposition} are satisfied.
\end{rem}

\section*{Acknowledgement}

I would like to thank Zolt\'an Szab\'o for his interest in this work
and Ko Honda and Tam\'as K\'alm\'an for inspiring conversations.

\bibliographystyle{amsplain}
\bibliography{topology}
\end{document}